\documentclass{fic-l}
\newtheorem{Thm}{Theorem}[section]
\newtheorem{Cor}[Thm]{Corollary}
\newtheorem{Prop}[Thm]{Proposition}
\newtheorem{Lemma}[Thm]{Lemma}
\theoremstyle{definition}
\newtheorem{Dfn}[Thm]{Definition}

\theoremstyle{remark}
\newtheorem{rem}[Thm]{Remark}

\def\ca{{\mathcal A}}
\def\cb{{\mathcal B}}

\def\cd{{\mathcal D}}

\def\ch{{\mathcal H}}
\def\cai{{\mathcal I}}

\def\ck{{\mathcal K}}
\def\cl{{\mathcal L}}
\def\cam{{\mathcal M}}

\def\bn{{\mathbb N}}
\def\br{{\mathbb R}}

\def\a{\alpha}
\def\b{\beta}
\def\g{\gamma}        
\def\d{\delta}        
\def\eps{\varepsilon}

\def\l{\lambda}       
\def\m{\mu}

\def\r{\rho}
\def\s{\sigma}
     \def\S{\Sigma}
\def\t{\tau}

\def\f{\varphi}
\def\o{\omega}        \def\O{\Omega}

\def\imply{\Rightarrow}

\def\ov{\overline}

\def\itm#1{\item[($#1$)]}

\DeclareMathOperator{\Lim}{Lim}

\DeclareMathOperator{\Tr}{Tr}
\DeclareMathOperator{\diam}{diam}
\DeclareMathOperator{\vol}{vol}

\def\ubd#1{\overline{d_{B}}(#1)} 
\def\lbd#1{\underline{d_{B}}(#1)}

\def\cpt{\ck(\ch)}

\def\subc{\underline{c}}
\def\supc{\overline{c}}
\def\subd{\underline{d}}
\def\supd{\overline{d}}

\def\oinf#1{{\rm ord}(#1)}

\begin{document}
\title{Fractals in Noncommutative Geometry}
 \author{Daniele Guido}
 \address{Dipartimento di Matematica,\\ Universit\`a della Basilicata,\\ 
 I--85100 Potenza, Italy.}
\email{guido@unibas.it}
 \author{Tommaso Isola}
 \address{Dipartimento di Matematica,\\ Universit\`a di Roma ``Tor Vergata'',\\ 
 I--00133 Roma, Italy.}
 \email{isola@mat.uniroma2.it}
\subjclass[2000]{}
\date{January 30, 2001}

 \begin{abstract}
	 To any spectral triple $(\ca,D,\ch)$ a dimension $d$ is associated, 
	 in analogy with the Hausdorff dimension for metric spaces. Indeed 
	 $d$ is the unique number, if any, such that $|D|^{-d}$ has non 
	 trivial logarithmic Dixmier trace. Moreover, when $d\in(0,\infty)$, 
	 there always exists a singular trace which is finite nonzero on 
	 $|D|^{-d}$, giving rise to a noncommutative integration on $\ca$.
	 
	 Such results are applied to fractals in $\br$, using Connes' 
	 spectral triple, and to limit fractals in $\br^{n}$, a class 
	 which generalises self-similar fractals, using a new spectral 
	 triple.  The noncommutative dimension or measure can be 
	 computed in some cases.  They are shown to coincide with the 
	 (classical) Hausdorff dimension and measure in the case of 
	 self-similar fractals.
 \end{abstract}
\maketitle

 \section{Introduction.}\label{sec:intro}
 
 This paper is both a survey and an announcement of results concerning 
 singular traces on $\cb(\ch)$, and their application to the study of 
 fractals in the framework of noncommutative geometry.
 
 Alain Connes' noncommutative geometry is a relatively young 
 discipline founded some twenty years ago, but it is rapidly 
 developing both in theory and the applications (see {\it e.g.} the 
 books by Connes \cite{Co}, Gracia-Bondia et al.  \cite{Varilly}, 
 Connes-Moscovici recent papers, \cite{CM}, etc.).  In all of 
 them, Dixmier logarithmic trace (or its companion Wodzicki 
 noncommutative residue \cite{Wo}) plays an important role, as 
 providing the proper analogue of integration in the noncommutative 
 context.
 
 One aspect of noncommutative geometry, or, more precisely, of the 
 notion of spectral triple, is that it is broad enough to treat also 
 commutative singular spaces, which are too irregular to be treated 
 with the instruments of Riemannian geometry, as {\it e.g.} fractals 
 (cf.  Connes' book).  We observe that, in this case, a lot of 
 different singular traces, more general than the logarithmic one, 
 naturally appear.  Indeed, in \cite{AGPS}, completing results 
 previously obtained by Dixmier \cite{Dixmier} and Varga \cite{Varga}, 
 singular traces on $\cb(\ch)$ have been studied and classified, and 
 have been discovered singular traces whose natural domain lies inside 
 $\cl^{1}$.
 
 Our approach to the study of fractals by means of noncommutative 
 geometry recovers the known results on Hausdorff dimension and 
 measure for the class of self-similar fractals.  Moreover it has 
 motivated the definition of dimension and of Hausdorff (and 
 Hausdorff-Besicovitch) measure in the abstract setting of spectral 
 triples, because of the strong analogies with the classical case.
 
 One virtue of working with general singular traces, and not just with 
 the logarithmic one, is that we can see the dimension of a spectral 
 triple as a number which produces a noncommutative measure, {\it 
 i.e.} a linear functional on the spectral triple.  Such functional is 
 not, in general, based on the logarithmic trace.  Moreover, these 
 more general singular traces appear naturally in the study of the 
 class of limit fractals.
 
 The paper is organised as follows.  The first section contains a 
 survey of known results on singular traces on $\cb(\ch)$ 
 \cite{AGPS,Dixmier,Varga}.  Moreover the notions of order and 
 exponent of singular traceability of a compact operator are 
 introduced, and some of their recently found properties are 
 described.  Examples illustrating the notion of exponent of singular 
 traceability constitute section two.
 
 In the third section, the notion of dimension of a spectral triple is 
 introduced, together with the associated Hausdorff-Besicovitch 
 functionals, which are constructed by means of singular traces.
 
 Then we recall known results on fractal sets, and introduce a new 
 class of fractals, which we call limit fractals.
 
 The fifth section contains an announcement of our results on 
 fractals in $\br$, approached by means of the spectral triple 
 introduced by Connes.
 
 In the last section, we associate a spectral triple to limit fractals in 
 $\br^{n}$, which is  different from Connes' one. In the self-similar 
 case, this triple recovers Hausdorff dimension and measure. Partial 
 results are described in the general case.

 
 \section{Singular traces on the compact operators of a Hilbert space.}

 In this section we present the theory of singular traces on 
 $\cb(\ch)$ as it was developed by Dixmier \cite{Dixmier}, who first 
 showed their existence, and then in \cite{Varga}, \cite{AGPS} and 
 \cite{GI5}.

 A singular trace on $\cb(\ch)$ is a tracial weight vanishing on the 
 finite rank projections.  Any tracial weight is finite on an ideal 
 contained in $\ck(\ch)$ and may be decomposed as a sum of a singular 
 trace and a multiple of the normal trace.  Therefore the study of 
 (non-normal) traces on $\cb(\ch)$ is the same as the study of 
 singular traces.  Moreover, making use of unitary invariance, a 
 singular trace should depend only on the eigenvalue asymptotics, 
 namely, if $A$ and $B$ are positive compact operators on $\ch$ and 
 $\m_n(A)=\m_n(B)+o(\m_n(B))$, $\m_n$ denoting the $n$-th eigenvalue, 
 then $\t(A)=\t(B)$ for any singular trace $\t$.  The main problem 
 about singular traces is therefore to detect which asymptotics may be 
 ``resummed'' by a suitable singular trace, that is to say, which 
 operators are singularly traceable.

 In order to state the most general result in this respect we need 
 some notation.  Let $A$ be a compact operator.  Then we denote by 
 $\{\m_n(A)\}$ the sequence of the eigenvalues of $|A|$, arranged in 
 non-increasing order and counted with multiplicity.  We consider also 
 the (integral) sequence $\{S_n(A)\}$ defined as follows:
 $$
 S_n(A):=
 \begin{cases}
	\sum_{k=1}^{n} \m_k(A) & A\notin\cl^1 \cr 
	\sum_{k=n+1}^{\infty}\m_k(A) & A\in\cl^1,
 \end{cases} 
 $$
 where $\cl^1$ denotes the ideal of trace-class operators.  We call a compact 
 operator {\it singularly traceable} if there exists a 
 singular trace which is finite non-zero on $|A|$.  We observe that 
 the domain of such singular trace should necessarily contain the 
 ideal $\cai(A)$ generated by $A$.  A compact operator is called {\it 
 eccentric} if
 \begin{equation}\label{eq:1.1}
	\frac{S_{2n_k}(A)}{S_{n_k}(A)}\to1
 \end{equation}
 for a suitable subsequence $n_k$.  Then the following theorem holds.

 \begin{Thm}\label{eccsingtrac}
	 A positive compact operator $A$ is singularly traceable $iff$ it 
	 is eccentric.  In this case there exists a sequence $n_k$ such 
	 that condition~(\ref{eq:1.1}) is satisfied and, for any 
	 generalised limit $\Lim_{\omega}$ on $\ell^\infty$, the positive 
	 functional
	 $$
	 \tau_\omega (B) = 
	 \begin{cases} 
		 \Lim_{\omega}\left(\left\{\frac{S_{n_k}(B)}{S_{n_k}(A)}\right\}\right) 
		 &\quad B \in \cai(A)_+ \\
		 +\infty&\quad B \not\in \cai(A),\ B>0,
	 \end{cases}
	 $$
	 is a singular trace whose domain is the ideal $\cai(A)$ generated 
	 by $A$. 
 \end{Thm}

 Now we give a sufficient condition to ensure eccentricity.  It is 
 based on the notion of order of infinitesimal.

 \begin{Dfn} 
	 \itm{i} For $A\in\cpt$ we define
	 {\it order of infinitesimal} of $A$ 
	 $$
	 \oinf{A} := \liminf_{n\to\infty} \frac{\log \m_{n}(A)}{\log 
	 (1/n)}\; ,
	 $$
	 \itm{ii}  $\a\in(0,\infty)$ is called an {\it exponent of singular traceability} 
	 for $A\in\cpt$ if there is a singular trace $\t$ on $\cb(\ch)$ such 
	 that $\t(A^{\a})=1$.
 \end{Dfn}
 
 \begin{rem}\label{rem1.6}
	 \itm{i} In \cite{GI5} we used ${\rm ord}_\infty$ instead of ${\rm ord}$.
	 \itm{ii} If $A\in\cpt_{+}$, then for any $\a>0$, $\oinf{A^{\a}}= \a\ 
	 \oinf{A}$.
 \end{rem}
   
 \begin{Thm}\label{inf-ecc} 
	 \itm{i} Let $A\in\cpt$ be s.t. $\oinf{A}=1$.  Then $A$ is 
	 eccentric.
	 \itm{ii} If $\oinf{A}\in(0,\infty)$, then $\oinf{A}^{-1}$ is an 
	 exponent of singular traceability. 
 \end{Thm}
 
 Because of its importance in determining the eccentricity property of 
 an operator, therefore the existence of a non-trivial singular trace, 
 we give alternative ways of computing the order of an operator. Recall
 
 \begin{Dfn}
 	 \begin{align*}
		 \cl^{1,\infty} & :=\{A\in\ck(\ch):\sum_{k=1}^{n}\m_k(A)=O(\log n)\}, \\
  		 \cl_{0}^{1,\infty} & :=\{A\in\ck(\ch):\sum_{k=1}^{n}\m_k(A)=o(\log n)\}.
     \end{align*}
 \end{Dfn}
 
 \begin{Thm}\label{equi-ord} 
	 Let $A\in\cpt_{+}$.  Then $\oinf{A}= \inf \{ \a>0 : A^{\a}\in \cl^{1} 
	 \} = \inf \{ \a>0 : A^{\a}\in \cl^{1,\infty} \} = \sup \{ \a>0 : 
	 A^{\a}\in \cl_{0}^{1,\infty} \}$.
 \end{Thm}
 
 Now we associate to any compact operator $A$ two numbers, which give 
 bounds for singular traceability.  We denote by $\m_A$ the locally 
 constant function defined by $\m(x)\equiv\m_A(x):=\m_n$ when 
 $x\in[n,n+1)$, $n\in\bn$, and by $f\equiv f_A$ the increasing, 
 diverging function determined by $f(t)=-\log\m(e^t)$.

 \begin{Dfn}
	 Let $A$ be a compact operator, $f\equiv f_A$ the increasing, 
	 diverging function defined before.  Then we set
	 \begin{align*}
		\subc(A) & = \left(\lim_{h\to\infty}\limsup_{t\to\infty} 
		\frac{f(t+h)-f(t)}{h}\right)^{-1}, \\
		\supc(A) & = \left(\lim_{h\to\infty}\liminf_{t\to\infty} 
		\frac{f(t+h)-f(t)}{h}\right)^{-1}.
	 \end{align*}
 \end{Dfn}
 
 \begin{Thm}\label{subcsupc}
	 Let $A$ be a compact operator.  Then the two limits above exist, 
	 and, if $\a$ is an exponent of singular traceability, then 
	 necessarily $\a\in[\subc(A),\supc(A)]$.  In particular 
	 $\subc(A)\leq\oinf{A}^{-1}\leq\supc(A)$.
 \end{Thm}
 
 The first result on singular traceability is due to Dixmier, who 
 showed in \cite{Dixmier} that $\frac{S_{2n}(A)}{S_n(A)}\to1$ is a 
 sufficient condition for singular traceability when $A\not\in\cl^1$.  
 Then Varga proved that the eccentricity condition is necessary and 
 sufficient when $A\not\in\cl^1$ \cite{Varga}.  Finally it was 
 observed in \cite{AGPSa} that singular traces may be non-trivial on 
 trace-class operators, while Theorem \ref{eccsingtrac} in the 
 previous form is contained in \cite{AGPS}.  Theorem \ref{inf-ecc} is 
 in \cite{GI5}, while the proof of Theorems \ref{equi-ord}, 
 \ref{subcsupc} will appear in \cite{GI9}.

 \medskip

 \section{Examples}
 
 This section is devoted to some examples, where the necessary 
 condition in Theorem \ref{subcsupc} is sufficient.  In particular in 
 the first class of examples $\subc$ and $\supc$ are finite non-zero, 
 and the exponents of singular traceability are exactly the elements 
 of $[\subc,\supc]$.  In the second class of examples $\subc=0$ and 
 $\supc=\infty$, and all positive numbers are exponents of singular 
 traceability.

 \subsection{On a class of operators for which all 
 $\g\in[\subc,\supc]$ are indices of singularly traceability.}

 In this subsection we will use the following notation:
 \begin{align*}
	 \s^{(\g)}(x)&:=\int_1^x     \m(y)^{\g}dy, \qquad \s(x) = \s^{(1)}(x)\cr
	 s^{(\g)}(x)&:=\int_x^\infty\m(y)^{\g}dy, \qquad s(x) = s^{(1)}(x)
 \end{align*}
 and the property (analogous to Theorem \ref{eccsingtrac}, cf.  
 \cite{GI1}) that $T$ is singularly traceable if and only if 1 is a 
 limit point of $\frac{\s(x)}{\s(2x)}$ or $\frac{s(x)}{s(2x)}$ as the 
 case may be.

 Let us choose two numbers $0<\b\leq\a$ and a non-decreasing sequence 
 $\{a_{n}\}_{n\in\bn}$, $a_{0}=0$, $a_{1}>0$.
 
 Set 
 \begin{align*}
 	b_{n}&=\sum_{k=0}^{n}a_{k}\cr
	\f(t)&=
	\begin{cases}
		\a\quad& t\in [b_{n},b_{n+1}), n $ even$\cr
		\b\quad& t\in [b_{n},b_{n+1}), n $ odd$
	\end{cases}
	\cr
	f(t)&=\int_{0}^{t}\f(s)ds
 \end{align*}
 
 Observe that $f$ is nondecreasing and goes to $\infty$ as 
 $t\to\infty$; as a consequence, $\m(x)=e^{-f(\log x)}$ is 
 nonincreasing and goes to 0.  We choose a compact operator $T$ such 
 that $\m_{n}(T)=\m(n)$, and prove the following.

 \begin{Thm}\label{exa1} 
	 $T^{\g}$ is singularly traceable {\it iff} $\g\in[\subc(T),\supc(T)]$
 \end{Thm}
 
 \begin{Lemma}\label{Lemma1} 
	If $\sup a_{n}=\lim a_{n}=a<\infty$, then $\subc=\supc$.
 \end{Lemma}
 \begin{proof}  
	If $k>1$, $\exists n_{0}$ such that $a\geq a_{n}>a(1-1/k)$ for all 
	$n\geq n_{0}$.  Then, on an interval $[t,t+ka]$, $t\geq n_{0}$, 
	there are $k/2\pm1$ intervals where $\f=\a$ and $k/2\pm1$ 
	intervals where $\f=\b$, hence
	\begin{align*}
		\frac{f(t+ka)-f(t)}{ka}
		&=\frac{1}{ka}\int_{t}^{t+ka}\f(s)ds\cr
		&=\frac{1}{ka}\left(a\a\left(\frac{k}{2}\pm 1\right)
		 +a\b\left(\frac{k}{2}\pm 1\right)\right)\cr
		&=\frac{\a+\b}{2}+\frac{\pm\a\pm\b}{k}
	\end{align*}
	which implies $\subc=\supc=\left(\frac{\a+\b}{2}\right)^{-1}$.
 \end{proof}

 \begin{Lemma} 
	 If $\sup a_{n}=\infty$ then $\subc=1/\a$, $\supc=1/\b$.
 \end{Lemma}
 \begin{proof}  
	For all $h>0$ $\exists n_{0}$ such that $n\geq n_{0}$ implies $ a_{n}>h$; 
	hence $b_{n+1}-b_{n}=a_{n+1}>h$ and 
	$$
	\frac{f(b_{n}+h)-f(b_{n})}{h}=
	\begin{cases}
		\a\quad& n $ even$\cr
		\b\quad& n $ odd$
	\end{cases}
	$$
	which implies $\supc^{-1}\leq\b$, $\subc^{-1}\geq\a$.  On the 
	other hand $\b\geq\frac{f(t+h)-f(t)}{h}\geq\a$, which implies the 
	thesis.
 \end{proof}

 \begin{proof} ({\it of Theorem \ref{exa1}}.)  If $a_n\to a<\infty$, 
 then $\subc=\supc=\oinf{T}^{-1}$ by Theorem \ref{subcsupc} and Lemma 
 \ref{Lemma1}, hence the thesis follows by Theorem \ref{inf-ecc}.  The 
 same argument applies if $\a=\b$, therefore we may assume $\a<\b$.  
 Because of Theorem \ref{subcsupc}, the claim is proved if we show 
 that, assuming $a_n\to\infty$ and $\a<\b$, $T^{\g}$ is singularly 
 traceable for any $\g\in[\subc,\supc]$.

 Observations:
 \begin{itemize}
	 \item[$(1)$] $t\in[b_{2n},b_{2n+1})$ implies
	 \begin{align*}	
		 f(t)&=f(b_{2n})+\int_{b_{2n}}^{t}\f=f(b_{2n})+(t-b_{2n})\a\cr 
		 &=f(b_{2n+1})-\int_{t}^{b_{2n+1}}\f=f(b_{2n+1})-(b_{2n+1}-t)\a.
	 \end{align*}
	 \item[$(2)$] $t\in[b_{2n-1},b_{2n})$ implies
	 \begin{align*}	
		 f(t)&=f(b_{2n-1})+(t-b_{2n-1})\b\cr 
		 &=f(b_{2n})-(b_{2n}-t)\b.
	 \end{align*}
 \end{itemize}

 Setting $x_{n}=e^{b_{n}}$,  we get
 $$
 \m(x)=
 \begin{cases}
	\m(x_{2n})\left(\frac{x_{2n}}{x}\right)^{\a}
	=\m(x_{2n+1})\left(\frac{x_{2n+1}}{x}\right)^{\a}
	& x\in[x_{2n},x_{2n+1}]\cr
	\m(x_{2n-1})\left(\frac{x_{2n-1}}{x}\right)^{\b}
	=\m(x_{2n})\left(\frac{x_{2n}}{x}\right)^{\b}
	& x\in[x_{2n-1},x_{2n}].
 \end{cases}
 $$
 Now assume $\subc<\g<\supc$, which is equivalent to $\a\g-1>0$ and 
 $\b\g-1<0$.  We shall show that, for any such $\g$, $1$ is a limit 
 point of ${\s^{(\g)}(\l x)}/{\s^{(\g)}(x)}$ and of 
 ${s^{(\g)}(x/\l)}/{s^{(\g)}(x)}$ for any $\l>1$ (but we only need 
 the case $\l=2$). Indeed
 $$
 0<\frac{\s^{(\g)}(\l x_{2n+1})}{\s^{(\g)}(x_{2n+1})}-1 \leq 
 \frac{\int_{x_{2n+1}}^{\l x_{2n+1}} 
 \left(\frac{x_{2n+1}}{x}\right)^{\b\g}dx} {\int_{x_{2n}}^{x_{2n+1}} 
 \left(\frac{x_{2n+1}}{x}\right)^{\a\g}dx} 
 =\frac{\a\g-1}{1-\b\g}\quad\frac{\l^{1-\b\g}-1}{e^{a_{2n+1}(\a\g-1)}-1}\to0
 $$

 $$
 0<\frac{s^{(\g)}(\frac{x_{2n-1}}{\l})}{s^{(\g)}(x_{2n-1})}-1 \leq 
 \frac{\int_{\frac{x_{2n-1}}{\l}}^{x_{2n-1}} 
 \left(\frac{x_{2n-1}}{x}\right)^{\a\g}dx} {\int_{x_{2n-1}}^{x_{2n}} 
 \left(\frac{x_{2n-1}}{x}\right)^{\b\g}dx } 
 =\frac{1-\b\g}{\a\g-1}\quad\frac{\l^{\a\g-1}-1}{e^{a_{2n}(1-\b\g)}-1}\to0
 $$

 Let now $\g=\subc=1/\a$; then

 $$
 0<\frac{\s^{(1/\a)}(\l x_{2n+1})}{\s^{(1/\a)}(x_{2n+1})}-1 \leq 
 \frac{\int_{x_{2n+1}}^{\l x_{2n+1}} 
 \left(\frac{x_{2n+1}}{x}\right)^{\b/\a}dx} {\int_{x_{2n}}^{x_{2n+1}} 
 \frac{x_{2n+1}}{x}dx} 
 =\frac{1}{1-\b/\a}\quad\frac{\l^{1-\b/\a}-1}{a_{2n+1}}\to0
 $$

 $$
 0<\frac{s^{(1/\a)}(\frac{x_{2n-1}}{\l})}{s^{(1/\a)}(x_{2n-1})}-1 \leq 
 \frac{\int_{\frac{x_{2n-1}}{\l}}^{x_{2n-1}} \frac{x_{2n-1}}{x}dx} 
 {\int_{x_{2n-1}}^{x_{2n}} \left(\frac{x_{2n-1}}{x}\right)^{\b/\a}dx } 
 =(1-\b/\a)\quad\frac{\log\l}{e^{a_{2n}(1-\b\g)}-1}\to0
 $$

 Finally, let $\g=\supc=1/\b$; then

 $$
 0<\frac{\s^{(1/\b)}(\l x_{2n+1})}{\s^{(1/\b)}(x_{2n+1})}-1 \leq 
 \frac{\int_{x_{2n+1}}^{\l x_{2n+1}} \frac{x_{2n+1}}{x}dx} 
 {\int_{x_{2n}}^{x_{2n+1}} \left(\frac{x_{2n+1}}{x}\right)^{\a/\b}dx} 
 =(\a/\b-1)\quad\frac{\log\l}{e^{a_{2n+1}(\a/\b-1)}-1}\to0
 $$

 $$
 0<\frac{s^{(1/\b)}(\frac{x_{2n-1}}{\l})}{s^{(1/\b)}(x_{2n-1})}-1 \leq 
 \frac{\int_{\frac{x_{2n-1}}{\l}}^{x_{2n-1}} 
 \left(\frac{x_{2n-1}}{x}\right)^{\a/\b}dx} {\int_{x_{2n-1}}^{x_{2n}} 
 \frac{x_{2n-1}}{x}dx} 
 =\frac{1}{\a/\b-1}\quad\frac{\l^{\a/\b-1}-1}{a_{2n}}\to0
 $$
 \end{proof}

 \begin{rem} 
	 It may happen that $\subc<\subd$ and $\supd<\supc$, where
	 $$
	 \subd=\left(\limsup_{t\to\infty}\frac{f(t)}{t}\right)^{-1},\qquad 
	 \supd=\left(\liminf_{t\to\infty}\frac{f(t)}{t}\right)^{-1}=\oinf{T}^{-1}.
	 $$

	 Choose $a_{n}=n$, $\b<\a$.  Then 
	 $b_{n}=\sum_{k=0}^{n}k=\frac{n(n+1)}{2}$.  If 
	 $t\in[b_{2n},b_{2n+1}]$, then
	 $$
	 f(b_{2n})\leq f(t)\leq f(b_{2n+2})=f(b_{2n})+\a(n+1)+\b(n+2),
	 $$
	 Hence
	 $$
	 \frac{f(b_{2n})}{(n+1)(2n+1)}\leq\frac{f(t)}{t} 
	 \leq\frac{f(b_{2n})}{n(2n+1)}+ \frac{\a(n+1)+\b(n+2)}{n(2n+1)}
	 $$
	 Finally, 
	 $$
	 f(b_{2n})=\sum_{j=1}^{n}\b(2j)+\sum_{j=1}^{n}\a(2j-1) 
	 =n(n+1)\b+n^{2}\a,
	 $$
	 which implies $\frac{f(b_{2n})}{2n^{2}}\to\frac{\a+\b}{2}$, 
	 therefore $\lim\frac{f(t)}{t}=\frac{\a+\b}{2}$.
 \end{rem}

 \subsection{On a class of operators for which all 
 positive numbers are indices of singularly traceability.}

 Choose an increasing sequence $b_n$ such that $e^{b_n}\in\bn$, 
 $b_{n+1}-b_{n}\to\infty$, $b_{0}=0$, and set
 $$
 f(t)=b_{n},\quad b_{n-1}<t\leq b_{n}.
 $$
 As before, set $\m(x)=e^{-f(\log x)}$, namely
 $$
 \m(x)=\frac{1}{x_{n}},\quad x_{n-1}<x\leq x_{n},
 $$
 where $x_{n}=e^{b_{n}}\in\bn$, hence $\frac{x_{n}}{x_{n+1}}\to0$.  We 
 choose a compact operator $T$ such that $\m_{n}(T)=\m(n)$, and prove 
 the following.

 \begin{Thm}\label{exa2} 
	 $T^{\a}$ is singularly traceable for any $\a>0$.
 \end{Thm}

 The proof of this statement requires some steps.  First we observe 
 that since $\forall h>0$ $\exists n_{0}$ such that $n>n_{0}\imply 
 b_{n+1}>b_{n}+h$, we have
 \begin{align*}
	 \frac{f(b_{n+1})-f(b_{n+1}-h)}{h}=0,\quad n>n_0 &\imply 
	 \liminf_{t\to\infty} \frac{f(t+h)-f(t)}{h}=0,\cr 
	 \frac{f(b_{n}+h)-f(b_{n})}{h}=\frac{b_{n+1}-b_{n}}{h},\quad n>n_0 
	 &\imply \limsup_{t\to\infty} \frac{f(t+h)-f(t)}{h}=+\infty,
 \end{align*}
 namely $\subc=0$, $\supc=\infty$. Also
 $$
 \oinf{T}=\liminf_{t\to\infty}\frac{f(t)}{t} 
 =\lim_{t\to\infty}\frac{f(b_{n})}{b_{n}}=1.
 $$

 \begin{Prop} 
	 Let $A$ be a compact operator.  If 
	 $\liminf\frac{\m_{n+1}(A)}{\m_{n}(A)}=0$, then $A^{\a}$ is 
	 singularly traceable for any $\a<(\oinf{A})^{-1}$.  
 \end{Prop}

 \begin{proof} 
	 By Theorem \ref{equi-ord}, when $\a<(\oinf{A})^{-1}$, $A^{\a}$ is 
	 not trace class; moreover $\oinf{A^{\a}}=\a\oinf{A}$ and 
	 $\m_{n}(A^{\a})=(\m_{n}(A))^{\a}$.  Therefore we may assume $A$ 
	 not to be trace class and $\a=1$.  Then, let $n_{k}$ be such that 
	 $\frac{\m_{n_{k}+1}(A)}{\m_{n_{k}}(A)}\to0$.  We have
	 $$
	 1\leq\frac{\s_{2n_{k}}}{\s_{n_{k}}} 
	 =1+\frac{\sum_{j=n_{k}+1}^{2n_{k}}\m_{j}}{\sum_{j=1}^{n_{k}}\m_{j}} 
	 \leq1+\frac{n_{k}\m_{n_{k}+1}}{n_{k}\m_{n_{k}}}\to1.
	 $$
	 The thesis follows by Theorem \ref{eccsingtrac}.
 \end{proof}

 \begin{Cor} 
	 $T^{\a}$ is singularly traceable for $\a\in(0,1)$.
 \end{Cor}
 \begin{proof} 
	 Indeed
	 $$
	 \frac{\m(x_{n+1})}{\m(x_{n})}=e^{-(b_{n+1}-b_{n})}\to0.
	 $$
 \end{proof}

 \begin{Lemma}
	 $\frac{s^{(\a)}(x_{n+1})}{x^{1-\a}_{n+1}}\to0$, for any $\a>1$.
 \end{Lemma}
 \begin{proof}
	 First we show that
	 \begin{equation}\label{eq1}
		 \lim_{n\to\infty}\frac{\sum_{k=n+1}^{\infty}x_{k}^{-\eps}}{x_{n}^{-\eps}}=0, 
		 \qquad \forall\eps>0.
	 \end{equation}
	 Indeed
	 $$
	 \frac{\sum_{k=n+1}^{\infty}x_{k}^{-\eps}} 
	 {\sum_{k=n}^{\infty}x_{k}^{-\eps}} =\frac 
	 {\sum_{k=n}^{\infty}\left(\frac{x_{k+1}}{x_{k}}\right)^{-\eps}x_{k}^{-\eps}} 
	 {\sum_{k=n}^{\infty}x_{k}^{-\eps}} \leq\sup_{k\geq 
	 n}\left(\frac{x_{k+1}}{x_{k}}\right)^{-\eps}\to0.
	 $$
	 Therefore
	 $$
	 \frac{x_{n}^{-\eps}}{\sum_{k=n+1}^{\infty}x_{k}^{-\eps}} 
	 =\frac{\sum_{k=n}^{\infty}x_{k}^{-\eps}} 
	 {\sum_{k=n+1}^{\infty}x_{k}^{-\eps}}-1\to\infty.
	 $$
	 Now we observe that
	 $$
	 s^{(\a)}(x_{n+1})=\sum_{k=x_{n+1}+1}^{\infty}\m_{k}^{\a} 
	 =\sum_{p=n+1}^{\infty}x_{p+1}^{(1-\a)}\left(1-\frac{x_{p}}{x_{p+1}}\right) 
	 \leq\sum_{p=n+1}^{\infty}x_{p+1}^{(1-\a)}.
	 $$
	 The thesis follows from equation (\ref{eq1}).
 \end{proof}

 \begin{Prop} 
	 $T^{\a}$ is singularly traceable for $\a>1$.
 \end{Prop}
 \begin{proof}
	 \begin{align*}
		 \frac{s_{2x_{n}}^{(\a)}}{s_{x_{n}}^{(\a)}} 
		 &=\frac{\sum_{k=2x_{n}+1}^{x_{n+1}}\m_{k}+s_{x_{n+1}}^{(\a)}} 
		 {\sum_{k=x_{n}+1}^{x_{n+1}}\m_{k}+s_{x_{n+1}}^{(\a)}}\cr 
		 &=\frac{(x_{n+1}-2x_{n})x_{n+1}^{-\a}+s_{x_{n+1}}^{(\a)}} 
		 {(x_{n+1}-x_{n})x_{n+1}^{-\a}+s_{x_{n+1}}^{(\a)}}\cr 
		 &=\frac{\left(1-2\frac{x_{n}}{x_{n+1}}\right) 
		 +\frac{s_{x_{n+1}}^{(\a)}}{x_{n+1}^{1-\a}}} 
		 {\left(1-\frac{x_{n}}{x_{n+1}}\right) 
		 +\frac{s_{x_{n+1}}^{(\a)}}{x_{n+1}^{1-\a}}}\to1.
	 \end{align*}
 \end{proof}

 \section{Some results on noncommutative geometric measure theory}\label{sec2}

 In this section we shall discuss a definition of dimension in 
 noncommutative geometry in the spirit of geometric measure theory.

 As it is known, the geometric measure for a noncommutative manifold 
 is defined via a singular trace applied to a suitable power of some 
 geometric operator (e.g. the Dirac operator of the spectral triple of 
 Alain Connes).  Connes showed that such procedure recovers the usual 
 volume in the case of compact Riemannian manifolds, and more 
 generally the Hausdorff measure in some interesting examples 
 \cite{Co}.

 Let us recall that $(\ca,D,\ch)$ is called a {\it spectral triple} 
 when $\ca$ is an algebra acting on the Hilbert space $\ch$, $D$ is a 
 self adjoint operator on the same Hilbert space such that $[D,a]$ is 
 bounded for any $a\in\ca$, and $D$ has compact resolvent.  In the 
 following we shall assume that $0$ is not an eigenvalue of $D$, the 
 general case being recovered by replacing $D$ with 
 $D|_{\ker(D)^\perp}$.  Such a triple is called $d^+$-summable, $d\in 
 (0,\infty)$, when $|D|^{-d}\in\cl^{1,\infty}$.

 The noncommutative version of the integral on functions is given by 
 the formula $a\mapsto\Tr_\omega(a|D|^{-d})$, where $\Tr_\omega$ is a 
 (logarithmic) Dixmier trace, i.e. a singular trace summing 
 logarithmic divergences.  Of course the preceding formula does not 
 guarantee the non-triviality of the integral, and in fact 
 cohomological assumptions in this direction have been considered 
 \cite{Co}.  We are interested in different conditions for 
 non-triviality.  In this connection, we observe that the previous 
 noncommutative integration is always trivial when $|D|^{-d}$ belongs 
 to $\cl^{1,\infty}_0$.

 \begin{Prop} 
	 Let $(\ca,D,\ch)$ be a spectral triple.  If $d$ is an exponent of 
	 singular traceability for $|D|^{-1}$, namely there is a singular 
	 trace $\t$ which is non-trivial on the ideal generated by 
	 $|D|^{-d}$, then the functional $a\mapsto\t(a 
	 |D|^{-d})$ is a non-trivial trace state on the algebra $\ca$. \\
	 We call it a Hausdorff-Besicovitch functional on $(\ca,D,\ch)$. Under 
	 suitable conditions (see \cite{CiGS1}) it gives rise to a trace on 
	 $\O\ca$.
 \end{Prop}

 \begin{rem}
	Any such trace is a candidate for a geometric measure in 
	noncommutative geometry.  Indeed, when $(\ca,D,\ch)$ is associated 
	to an $n$-dimensional compact manifold $M$, or to the fractal sets 
	in \cite{Co}, the singular trace is the logarithmic Dixmier trace, 
	and the associated functional corresponds to the Hausdorff 
	measure.  Therefore the following definition is natural.
 \end{rem}
 
  \begin{Dfn} 
	  Let $(\ca,D,\ch)$ be a spectral triple, $\Tr_{\o}$ the 
	  logarithmic Dixmier trace.  \itm{i} We call $\a$-dimensional 
	  Hausdorff functional the map $a\mapsto Tr_\omega(a |D|^{-\a})$; 
	  \itm{ii} we call (Hausdorff) dimension of the spectral triple 
	  the number
	  $$
	  d(\ca,D,\ch) = \inf \{ d>0: |D|^{-d} \in \cl^{1,\infty}_{0} \} = 
	  \sup \{ d>0: |D|^{-d} \not\in \cl^{1,\infty} \}.
	  $$
 \end{Dfn}
 
 \begin{Thm}\label{unisingtrac} 
	 Assume $d:=d(\ca,D,\ch)\in(0,\infty)$. Then 
	 \itm{i} $d$ is the unique exponent, if any, such that $\ch_{d}$ 
	 is non-trivial; 
	 \itm{ii} $d = \oinf{D^{-1}}^{-1} = \sup\{\g>0:|D|^{-\g} 
	 \notin\cl^1\}$; as a consequence it is an exponent of singular 
	 traceability;  
	 \itm{iii} let $\zeta(s)=\Tr(|D|^{-s})$, $s>d$.  If 
	 $(s-d)\zeta(s)\to L\in\br$ as $s\to d^+$, then the 
	 Hausdorff-Besicovitch functional associated with $|D|^{-d}$ is 
	 indeed the Hausdorff functional (up to the multiplicative 
	 constant $L$).
 \end{Thm}
 \begin{proof}
	 $(i)$ is in \cite{GI5}, $(ii)$ follows from Theorems 
	 \ref{inf-ecc} and \ref{equi-ord}, $(iii)$ follows by the 
	 Hardy-Littlewood Theorem \cite{Hardy}, cf.  Proposition 4, p.  
	 306, \cite{Co}.
 \end{proof}

 Let us observe that the $\a$-dimensional Hausdorff functional depends 
 on the generalised limit procedure $\omega$, however all such 
 functionals coincide on the elements $a$ of $\ca$ such that $a|D|^{-d}$ 
 is a measurable operator in the sense of Connes \cite{Co}.  As in the 
 commutative case, the dimension is the supremum of the $\a$'s such 
 that the $\a$-dimensional Hausdorff measure is everywhere infinite 
 and the infimum of the $\a$'s such that the $\a$-dimensional 
 Hausdorff measure is identically zero.

 Concerning the non-triviality of the $d$-dimensional Hausdorff 
 functional, we have the same situation as in the classical case.  
 Indeed, according to the previous result, a non-trivial Hausdorff 
 functional is unique but does not necessarily exist.  In fact, if the 
 eigenvalue asymptotics of $D$ is e.g. $n\log n$, the Hausdorff 
 dimension is one, but the 1-dimensional Hausdorff measure gives the 
 null functional.

 However, if we consider all singular traces, not only the logarithmic 
 ones, and the corresponding functionals on $\ca$, as we said, there 
 exists a non trivial functional associated with such a dimension, but 
 such property does not characterize this dimension, in general, namely the
 exponent of singular traceability is not necessarily unique, cf. the 
 examples in Section 2.
 
  \begin{Prop}\label{Prop:unique}
	  If $\subc(|D|^{-1})=\supc(|D|^{-1})\in(0,\infty)$, then 
	  $d(\ca,D,\ch)=\subc=\supc$ is the unique exponent of singular 
	  traceability of $D^{-1}$.  This is the case, in particular, if 
	  there exists 
	  $\lim\frac{\m_n(D^{-1})}{\m_{2n}(D^{-1})}\in(1,\infty)$.
  \end{Prop}
  \begin{proof} 
	  The first statement follows from \ref{subcsupc}.  Since 
	  $\lim\frac{\m_n(D^{-1})}{\m_{2n}(D^{-1})}=2^{\subc}=2^{\supc}$ 
	  if it exists, the second statement follows from the first, 
	  however it has been proved directly in \cite{GI5}.
  \end{proof}

 \begin{rem} For the spectral triples whose Dirac operator has a spectral 
	 asymptotics like $n^{\a}(\log n)^{\b}$, we have 
	 $d(\ca,D,\ch)=1/\a$, and the uniqueness result of Proposition 
	 \ref{Prop:unique} applies.  However, the nontrivial singular trace 
	 associated with $|D|^{-1/\a}$ by Theorem \ref{inf-ecc} is a 
	 logarithmic trace if and only if $\b=1$.  In this sense, the 
	 singular traces associated with a generic eccentric operator 
	 generalize the logarithmic Dixmier trace in the same way in which the 
	 Besicovitch measure theory generalizes the Hausdorff measure 
	 theory.
 \end{rem}
 
 \begin{rem}
	Contrary to the classical case, where there are sets with 
	non-trivial Hausdorff dimension but no non-trivial geometric 
	($i.e.$ Hausdorff or Hausdorff-Besicovitch) measure \cite{Rogers}, 
	p.  73, in the noncommutative context, if $d(\ca, D, \ch) \in 
	(0,\infty)$, there is always a non-trivial geometric measure, 
	whether Hausdorff or (the more general) Hausdorff-Besicovitch.
 \end{rem}

 \section{Fractals in $\br^{n}$. Classical aspects}\label{SecClassic}
 
 Let $(X,\r)$ be a metric space, and let $h:[0,\infty) \to [0,\infty)$ 
 be non-decreasing and right-continuous, with $h(0)=0$.  When 
 $E\subset X$, define, for any $\d>0$, $\ch^{h}_{\d}(E) := \inf \{ 
 \sum_{i=1}^{\infty} h(\diam A_{i}) : \cup_{i} A_{i} \supset E, \diam 
 A_{i} \leq \d \}$.  Then the {\it Hausdorff-Besicovitch (outer) 
 measure} of $E$ is defined as
 $$
 \ch^{h} (E) := \lim_{\d\to0}\ch^{h}_{\d}(E).
 $$
 If $h(t) = t^{\a}$, $\ch^{\a}$ is called {\it Hausdorff (outer) 
 measure} of order $\a>0$.
 
 The number 
 $$
 d_{H}(E) := \sup \{ \a>0 : \ch^{\a}(E) = +\infty \} = \inf \{ \a>0 : 
 \ch^{\a}(E) = 0 \}
 $$ 
 is called {\it Hausdorff dimension} of $E$.
 
 Let $N_{\eps}(E)$ be the least number of closed balls of radius 
 $\eps>0$ necessary to cover $E$.  Then the numbers
 $$
 \ubd{E} := \limsup_{\eps\to0^{+}} \frac{\log N_{\eps}(E)}{-\log 
 \eps}, \quad \lbd{E} := \liminf_{\eps\to0^{+}} \frac{\log 
 N_{\eps}(E)}{-\log \eps}
 $$ 
 are called upper and lower {\it box dimensions} of $E$.
 
 In case $X=\br^{N}$, setting $S_{\eps}(E) := \{ x\in\br^{N} : \r(x,E) 
 \leq \eps \}$, it is known that $\ubd{E} = N - \liminf_{\eps\to0^{+}} 
 \frac{\log \vol S_{\eps}(E)}{\log \eps}$ and $\lbd{E} = N - 
 \limsup_{\eps\to0^{+}} \frac{\log \vol S_{\eps}(E)}{\log \eps}$.  $E$ 
 is said {\it Minkowski measurable} if
 $$
 \lim_{\eps\to0^{+}} \frac{\log \vol S_{\eps}(E)}{\log \eps} = N-d 
 \quad\text{and}\quad \cam_{d}(E):=\lim_{\eps\to0^{+}} \frac{\vol 
 S_{\eps}(E)}{\eps^{N-d}}\in (0,\infty).
 $$ 
 $\cam_{d}(E)$ is called {\it Minkowski content} of $E$.

 \subsection{Selfsimilar fractals.}
 
 Let $\{w_{j} \}_{j=1,\ldots,p}$ be contracting similarities of 
 $\br^{N}$, $i.e.$ there are $\l_{j}\in(0,1)$ such that $\| w_{j}(x) - 
 w_{j}(y) \| = \l_{j} \|x-y\|$, $x,y\in \br^{N}$.  Denote by 
 $\ck(\br^{N})$ the family of all non-empty compact subsets of 
 $\br^{N}$, endowed with the Hausdorff metric, which turns it into a 
 complete metric space.  Then $W: K\in\ck(\br^{N}) \to \cup_{j=1}^{p} 
 w_{j}(K) \in \ck(\br^{N})$ is a contraction.
 \begin{Dfn}
	 The unique non-empty compact subset $F$ of $\br^{N}$ such that
	 $$
	 F = W(F) = \bigcup_{j=1}^{p} w_{j}(F)
	 $$ 
	 is called the {\it self-similar fractal} defined by $\{w_{j} 
	 \}_{j=1,\ldots,p}$.
 \end{Dfn}

 If we denote by $Prob_\ck(\br^{N})$ the set of probability measures 
 on $\br^{N}$ with compact support endowed with the Hutchinson metric, 
 $i.e.$ $d(\mu,\nu) := \sup \{ |\int fd\mu - \int fd\nu| : \|f\|_{Lip} 
 \leq 1 \}$, then the map
 $$
 \begin{matrix}
	 T : & Prob_\ck(\br^{N}) &\to &Prob_\ck(\br^{N})\cr 
	 &\mu&\mapsto&\sum_{j=1}^{p} \l_{j}^{s}\mu\circ w_{j}^{-1}
 \end{matrix} 
 $$
 is a contraction, where $s>0$ is the unique real number, called similarity
 dimension, satisfying $\sum_{j=1}^{p} \l_{j}^{s} =1$. We then observe that if
 $\m$ has support $K$, then $T\m$ has support $W(K)$. Since the sequence
 $W^n(K)$ is convergent, it turns out that it is bounded, namely there exists a
 compact set $K_0$ containing the supports of all the measures $T^n\m$. But on
 the space $Prob(K_0)$ the Hutchinson metric induces the weak$^*$ topology, and
 this space is compact in such topology, hence complete in the Hutchinson
 metric. Therefore there exists a fixed point of $T$ in $Prob_\ck(\br^{N})$,
 which is of course unique.

 \noindent{\bf{Open Set Condition.}} The similarities $\{w_{j} 
 \}_{j=1,\ldots,p}$ are said to satisfy the open set condition if 
 there is a non-empty bounded open set $V\subset \br^{N}$ such that 
 $\cup_{j=1}^{p} w_{j}(V) \subset V$ and $w_{i}(V)\cap w_{j}(V) 
 =\emptyset$, $i\neq j$.  In this case $d_{H}(F) = \lbd{F} = \ubd{F} = 
 s$, and the Hausdorff measure $\ch^{s}$ is non-trivial on $F$.  
 Therefore $\ch^{s}|_{F}$ is the unique (up to a constant factor) 
 Borel measure $\mu$, with compact support, such that $\mu(A) = 
 \sum_{j=1}^{p} \l_{j}^{s} \mu(w_{j}^{-1}(A))$, for any Borel subset 
 $A$ of $\br^{N}$.

 It has recently been proved \cite{Gatzuras} that, if the similarities 
 $\{w_{j} \}_{j=1,\ldots,p}$ satisfy the open set condition and $\log 
 \l_{1}$, \ldots $\log \l_{p}$ generate $(\br,+)$ as a minimal closed 
 subgroup, then $F$ is Minkowski measurable.
 
 \subsection{Limit fractals.}
 
 Several generalisations of the class of self-similar fractals have 
 been studied.  Here we propose a new one, that we call the class of 
 limit fractals.  For its construction we need the following theorem.
 
 \begin{Thm}\label{gencontraction}
 	Let $(X,\r)$ be a complete metric space, $T_{n}:X\to X$ be such that 
 	there are $\l_{n}\in(0,1)$ for which $\r(T_{n}x,T_{n}y) \leq \l_{n} 
 	\r(x,y)$, for $x,y\in X$. Assume $\sum_{n=1}^{\infty} 
 	\prod_{j=1}^{n} \l_{j} <\infty$, and there is $x\in X$ such that 
 	$\sup_{n\in\bn} \r(T_{n}x,x) <\infty$. Then
	\itm{i} $\sup_{n\in\bn} \r(T_{n}y,y) <\infty$, for any $y\in X$, 
	\itm{ii}	$\lim_{n\to\infty} T_{1}\circ T_{2}\circ\cdots\circ 
	T_{n}x=x_{0}\in X$ for any $x\in X$.
 \end{Thm}
 \begin{proof}
 	$(i)$ $\r(T_{n}y,y) \leq \r(T_{n}y,T_{n}x) + \r(T_{n}x,x) + \r(x,y) 
 	\leq (1+\l_{n})\r(x,y) + \r(T_{n}x,x)$, so that $\sup_{n\in\bn} 
 	\r(T_{n}y,y) \leq 2\r(x,y) + \sup_{n\in\bn} \r(T_{n}x,x) <\infty$.\\
	$(ii)$ Set $M:= \sup_{n\in\bn} \r(T_{n}x,x) <\infty$, and $S_{n} 
	:= T_{1}\circ T_{2}\circ\cdots\circ T_{n}$, $n\in\bn$.  As 
	$\r(S_{n+1}x, S_{n}x) \leq \l_{1}\l_{2}\cdots\l_{n}\r(T_{n+1}x,x) 
	\leq M \l_{1}\l_{2}\cdots\l_{n}$, there follows, for any 
	$p\in\bn$, $\r(S_{n+p}x,S_{n}x)\leq \r(S_{n+p}x,S_{n+p-1}x) + 
	\ldots + \r(S_{n+1}x,S_{n}x) \leq M \sum_{k=n}^{n+p-1}\prod_{j=1}^{k} 
	\l_{k} \leq M \sum_{k=n}^{\infty}\prod_{j=1}^{k}\l_{k} \to 0$, as 
	$n\to\infty$, that is $\{ S_{n}x \}$ is Cauchy in $X$. Therefore 
	there is $x_{0}\in X$ such that $S_{n}x\to x_{0}$. \\
	Let us prove that $x_{0}$ is independent of $x$. Indeed, if $y\in X$, 
	then $\r(S_{n}x,S_{n}y) \leq \l_{1}\l_{2}\cdots \l_{n} \r(x,y) \to 
	0$, as $n\to\infty$, so that $S_{n}x$ and $S_{n}y$ have the same 
	limit.
 \end{proof}
 
 \begin{rem} 
	 A sufficient condition for $\sum_{n=1}^{\infty} \prod_{j=1}^{n} 
	 \l_{j} <\infty$ to hold is $$\sup_{n\in\bn} \l_{n}<1.$$
 \end{rem}
 
 We now describe the class of limit fractals.  Let $\{w_{nj} \}$, 
 $n\in\bn$, $j=1,\ldots,p_{n}$, be contracting similarities of 
 $\br^{N}$, with contraction parameter $\l_{nj}\in(0,1)$.  Set $\Sigma 
 := \cup_{n\in\bn} \{\s:\{1,\ldots,n\}\to \bn : \s(k) \in 
 \{1,\ldots,p_{k}\}, k=1,\ldots,n \}$, and write $w_{\s} := 
 w_{1\s(1)}\circ w_{2\s(2)} \circ \cdots \circ w_{n\s(n)}$, for any 
 $\s\in\Sigma$, $|\s|=n$.  Assume $\sup_{n,j} \l_{nj} < 1$ and $\{ 
 w_{\s}(x) : \s\in\Sigma \}$ is bounded, for some (hence any) 
 $x\in\br^{N}$.  Then, by Theorem \ref{gencontraction}, the sequence of 
 maps $W_{n} : K\in\ck(\br^{N}) \to \cup_{j=1}^{p_{n}} w_{nj}(K) \in 
 \ck(\br^{N})$ is such that $\{ W_{1}\circ W_{2}\circ\cdots\circ 
 W_{n}(K) \}$ has a limit in $\ck(\br^{N})$, which is independent of 
 $K\in\ck(\br^{N})$.
 
 \begin{Dfn} 
	 The unique compact set $F$ which is the limit of $\{ W_{1}\circ 
	 W_{2}\circ\cdots\circ W_{n}(K) \}_{n\in\bn}$ is called the limit 
	 fractal defined by $\{w_{nj} \}$.  In the particular case that 
	 $\l_{nj}=\l_{n}$, $j=1,\ldots,p_{n}$, $n\in\bn$, $F$ is called a 
	 translation (limit) fractal.  The limit fractal $F$ is said to 
	 satisfy the {\it countably ramified open set condition} if there 
	 exists a nonempty bounded open set $V$ in $\br^n$ for which 
	 $w_{nj}(V) \subset V$ and $\ov{w_{ni}(V)}\cap \ov{w_{nj}(V)}$ is 
	 at most countable, for any $n$, $i\neq j$.
 \end{Dfn}
 
 As before, we may consider the action of the similarities on 
 measures, besides that on sets.  Given $s>0$ we set
 $$
 \begin{matrix}
	 T_n : & Prob_\ck(\br^{N}) &\to &Prob_\ck(\br^{N})\cr 
	 &\mu&\mapsto&\frac{1}{\sum_{j=1}^{p} \l_{nj}^{s}} \sum_{j=1}^{p} 
	 \l_{nj}^{s}\mu\circ w_{nj}^{-1}
 \end{matrix}
 $$
 and consider the sequence $\{T_{1}\circ T_{2}\circ\cdots\circ 
 T_{n}\m\}_{n\in\bn}$.  As before the supports of all such measures 
 are contained in a common compact set, therefore Theorem 
 \ref{gencontraction} applies and we get a unique limit measure 
 $\m_s$, depending on the chosen $s$.
 
 If the countably ramified open set condition holds, the sets 
 $w_{\s\cdot i}\ov{V}$, $w_{\s\cdot j}\ov{V}$ are essentially disjoint 
 when $i\ne j$, where $\s\cdot i$ is the concatenation of strings, and 
 are related by the similarity $w_{ni}\circ w_{nj}^{-1}$, $n=|\s|+1$, 
 therefore the measure $\m_s$ is the unique probability measure with 
 support the limit fractal $F$ which is homogeneous with parameter 
 $s$.  If $F$ is a translation fractal, $w_{ni}\circ w_{nj}^{-1}$ is 
 indeed an isometry, hence $\m_s$ is independent of $s$ and is the 
 unique  probability measure on $F$ which is  invariant under the 
 mentioned isometries.
 
 \subsection{Fractals in $\br$.}\label{4.3}
 
 We now specialise to subsets of $\br$ and survey some known results 
 on compact, totally disconnected subsets of $\br$, without isolated 
 points.  Let $F$ be such a set, and denote by $[a,b]$ the least 
 closed interval containing $F$.  Then $[a,b]\setminus F$ is the 
 disjoint union of open intervals $(a_{n},b_{n})$, where $b_{n}-a_{n} 
 \leq b_{n-1}-a_{n-1}$, $n\geq2$.  If $F$ has Lebesgue measure zero, 
 $i.e.$ if $\sum_{n=1}^{\infty} (b_{n}-a_{n}) = b-a$, then
 \begin{equation}\label{dimlin}
	 \ubd{F}=\limsup_{n\to\infty} \frac{\log n}{|\log(b_{n}-a_{n})|}.
 \end{equation}
 
 It has been proved in \cite{KaSa} that, when $F$ is a translation
 fractal in $\br$, there is a gauge function $h$ such that the
 corresponding Hausdorff-Besicovitch measure $\ch^{h}$ is non-trivial on $F$.
 Moreover, if $\lim_{t\to0} \frac{\log h(t)}{\log t} = \a$, then $d_{H}(F) =
 \a$. As a consequence, if the countably ramified open set condition holds,
 $\ch^{h}|_F$ coincides (up to a constant) with the limit measure $\m$ of the
 previous subsection.

 \section{Fractals in $\br$. Noncommutative aspects.}
 
 Let $F$ be a compact, totally disconnected subset of $\br$, without 
 isolated points, and let $a,b,a_n,b_n$ be as in subsection \ref{4.3}.
 Set $\ch_{n}:= \ell^{2}(\{a_{n},b_{n}\})$, $\ch := 
 \oplus_{n=1}^{\infty} \ch_{n}$,
 $$
 D_{n} := \frac{1}{b_{n}-a_{n}} 
 \begin{pmatrix}
 	0 & 1\\
	1 & 0
 \end{pmatrix}
 $$
 $D:= \oplus_{n=1}^{\infty} D_{n}$.  Consider the action of $C(F)$ on 
 $\ch$ by left multiplication: $(f\xi)(x)=f(x)\xi(x)$, $x\in 
 \cd:=\{a_n,b_n : n\in\bn\}$, and define $\ca := Lip(F)$.  Then
 
 \begin{Thm}\label{linfrac}\cite{Co}
 	\itm{i} $(\ca,D,\ch)$ is a spectral triple
	\itm{ii} the characteristic values of $D^{-1}$ are the numbers 
	$b_{n}-a_{n}$, $n\in\bn$, each with multiplicity 2.
 
 	If $F$ is Minkowski measurable, and has box dimension $d\in(0,1]$, 
 	then
	\itm{iii} $|D|^{-d}\in \cl^{1,\infty}$
	\itm{iv} $\Tr_{\o}(|D|^{-d}) = 2^{d}(1-d) \cam_{d}(F)$.
 \end{Thm}

 Statements $(iii)$ and $(iv)$ follow from results of Lapidus and Pomerance, 
 \cite{Lapidus}.
 \medskip
 
 Even if $F$ is not Minkowski measurable, we have
 
 \begin{Thm}
	$d(\ca,D,\ch) = \ubd{F}$.  Therefore, if $\ubd{F}\ne0$, we get a 
	Hausdorff-Besicovitch functional on the spectral triple, giving 
	rise to a non-trivial measure $\mu$ on $F$.
 \end{Thm}
 \begin{proof}
	 Follows by equation (\ref{dimlin}), Theorem \ref{unisingtrac} and 
	 Theorem \ref{linfrac}.(ii).
 \end{proof}
 
 \begin{rem}
 	If $F$ is a limit fractal with countably ramified open set 
 	condition, then $\mu$ can be explicitely computed, in particular it 
 	only depends on $F$. If $F$ is a translation limit fractal, then 
 	$\mu$ coincides with the measure in subsection 4.3.
 \end{rem}
 
 \begin{Thm}
 	Let $F$ be a self-similar fractal, and $s\in[0,1]$ its Hausdorff 
 	dimension. Then $s$ is the unique exponent of singular 
 	traceability for $D^{-1}$, and the Hausdorff-Besicovitch functional 
 	on the spectral triple corresponds to the $s$-dimensional Hausdorff 
 	measure on $F$.
 \end{Thm}
 \begin{proof}
	  Define
	  $$
	  S_{j}\xi (b) :=
	  \begin{cases}
		  \xi(w_{j}^{-1}(b)) & b\in w_{j}\cd \\
		  0 & b\not\in w_{j}\cd.
	  \end{cases}
	  $$
	  Then $S_{j}$ is an isometry and $|D|^{-d} = \sum_{j=1}^{p} 
	  \l_{j}^{d} S_{j} |D|^{-d} S_{j}^{*}$.  Therefore, if $d$ is an 
	  exponent of singular traceability for $|D|^{-1}$, the 
	  corresponding Hausdorff-Besicovitch functional is homogeneous of 
	  order $d$.  This implies that $d$ coincides with $s$, namely $s$ 
	  is the unique exponent of singular traceability, and the 
	  Hausdorff-Besicovitch functional corresponds to the 
	  $d$-dimensional Hausdorff measure.
  \end{proof}
 
  \begin{rem}
  	One can show that the noncommutative $s$-dimensional Hausdorff functional 
  	is nontrivial on $C(F)$. 
  \end{rem}

 \section{Fractals in $\br^n$}

 In this Section we treat the case of self-similar and limit fractals 
 in $\br^{n}$.  Here the construction of the spectral triple due to 
 Connes does not apply, and we propose a new construction for the 
 Dirac operator.  This construction is similar to Connes' in that it 
 is based on a discrete approximation of the fractal.  On the other 
 hand it differs from that of Connes since the eigenvalues are not 
 proportional to the size of the ``holes'' of the fractal, but to the 
 size of the remaining parts.  So in general the two constructions do 
 not agree on fractals in $\br$, even though we shall show in many 
 cases that they give rise to the same measure on the fractal.

 \subsection{Self-similar fractals.}

 Assume that $F$ is a self-similar fractal, constructed via the 
 similarities $w_{n}$, $n=1,\dots p$, satisfying open set condition 
 w.r.t the bounded open set $V$.  Choose two points $x,y\in V$ 
 and consider the points $x_{\s}:=w_{\s}x$, $y_{\s}:=w_{\s}y$, 
 $\s\in\S$.  Define the space 
 $\ch_{\s}:=\ell^{2}(\{x_{\s},y_{\s}\})$ with the operator
 $$
 D_{\s}:=\frac{1}{d(x_{\s},y_{\s})}
 \begin{pmatrix}
	0&1\\
	1&0
 \end{pmatrix}
 $$
 and set $\ch:=\oplus_{\s\in\S}\ch_{\s}$, $D:=\oplus_{\s\in\S}D_{\s}$.  
 Finally introduce the algebra $\ca$ of Lipschitz functions on $V$, 
 acting by left multiplication on $\ch$: 
 $$
 (f\xi)(x)=f(x)\xi(x),\quad x\in\{x_\s,y_\s:\s\in\S\}.
 $$
 In this way 
 $$
 [D,f]=\bigoplus_{\s\in\S}\frac{f(x_{\s})-f(y_{\s})}{d(x_{\s},y_{\s})}
 \begin{pmatrix}
	0&1\\
	-1&0
 \end{pmatrix},
 $$
 namely $(\ca,D,\ch)$ is a spectral triple.

 The following theorem holds:

 \begin{Thm}\label{ssRn}
	Let $(\ca,D,\ch)$ be the spectral triple associated with a 
	self-similar fractal $F$ with open set condition as above.  Then 
	the dimension $d$ of the triple coincides with the Hausdorff 
	dimension of $F$, and the noncommutative $d$-dimensional Hausdorff 
	functional corresponds to the classical Haudorff measure.
 \end{Thm}
 \begin{proof}
	 The proof will appear in \cite{GI9}
 \end{proof}

 \begin{rem}
	We note that Theorem \ref{ssRn} implies that the noncommutative 
	dimension and measure do not depend on the starting points $x,y$. 
	Moreover one can replace the pair with any finite family of pairs 
	without affecting the result.
	
	For translation self-similar fractals in $\br$, namely fractals where all 
	similarity parameters coincide, the starting pairs can be chosen in 
	such a way that the spectral triple of Connes coincides with ours.
	
	In general however this is not the case, and the distance induced by 
	our Dirac operator is different from the original one.
 \end{rem}

 \subsection{Limit fractals.}

 Let $F$ be a limit fractal as in Section \ref{SecClassic}.  The 
 spectral triple can be defined exactly as for the self-similar case.  
 Here we shall assume, besides open set condition, also countable 
 ramification, namely $w_{\s\cdot i}\overline{V}\cap w_{\s\cdot 
 j}\overline{V}$ is at most countable when $|\s|=n-1$, 
 $i,j=1,\dots,p_{n}$, $i\ne j$.

 Then the following holds:

 \begin{Thm}
	Let $(\ca,D,\ch)$ be the spectral triple associated with a limit 
	fractal $F$ with countably ramified open set condition as 
	above. Then for any exponent $\a$ of singular traceability for 
	$|D|^{-1}$ the corresponding noncommutative Hausdorff-Besicovitch 
	measure coincides with the limit measure on $F$ with scaling 
	parameter $\a$.
 \end{Thm}
 \begin{proof}
	 The proof will appear in \cite{GI9}
 \end{proof}

 Now we restrict to the class of translation fractals, namely limit 
 fractals for which $\l_{n,i}=\l_{n}$.  In this case we have a formula 
 for the spectral dimension:

 \begin{Thm}
	Let $(\ca,D,\ch)$ be the spectral triple associated with a 
	translation fractal $F$ with countably ramified open set 
	condition, where the similarities $w_{n,i}$, $i=1,\dots,p_{n}$ 
	have scaling parameter $\l_{n}$.  Then the spectral dimension is 
	given by the formula
	$$
	d=\limsup_{n}\frac{\sum_1^n\log p_k}{\sum_1^n\log 1/\l_k}.
	$$
	Moreover the measure corresponding to the associated 
	Hausdorff-Besicovitch functional is the unique probability measure 
	on $F$ invariant under the internal isometries of $F$.
 \end{Thm}
 \begin{proof}
	 The proof will appear in \cite{GI9}
 \end{proof}


\end{document}